\documentclass[letterpaper, 10 pt, conference]{ieeeconf}

\IEEEoverridecommandlockouts

\let\proof\relax

\usepackage{amsthm} 
\usepackage{graphicx}
\graphicspath{{./Figures/}} 
\usepackage{amsfonts,latexsym}
\usepackage{amsmath} 
\usepackage{amssymb} 
\usepackage[ruled, lined, ,longend, linesnumbered]{algorithm2e}
\SetKwInOut{Require}{Require} 
\usepackage{booktabs} 
\usepackage{threeparttable} 
\usepackage{tablefootnote} 
\usepackage{array}
\usepackage[overload]{empheq}
\usepackage{epstopdf} 

\def\epi#1{{\rm epi}(#1)}

\def\bx{\bar{x}}

\def\bm{\bar{m}}

\def\arg{{\rm{arg}}}

\def\dhh{{\nabla h}}
\def\dom#1{{\rm dom}(#1)}

\def\ri#1{{\rm ri}(#1)}

\newcommand {\bsis} {\left\{ \begin{array} }
	\newcommand {\esis} {\end{array}\right.}

\def\Sum#1#2{\sum\limits_{#1}^{#2}}
\def\set#1#2{\{ \; #1 \;:\;#2\;\}} 
\def\inR#1{\in\mathbb{R}^{#1}} 
\def\R{ {\rm \,I\!R} }    

\def\sp#1#2{\langle #1,#2\rangle }

\def\fracg#1#2{{\displaystyle{\frac{#1}{#2}}}}  

\newcommand{\cX}{{\mathcal{X}}}

\newtheorem{assumption}{Assumption}
\newtheorem{definition}{Definition}
\newtheorem{property}{Property}

\newcommand{\blista}{\renewcommand{\labelenumi}{(\roman{enumi})} 
\begin{enumerate}}
\newcommand{\elista}{\end{enumerate} \renewcommand{\labelenumi}{\arabic{enumi}.}}
\newcommand {\T}{^{\top}} 

\makeatletter 
\newcommand{\thickhline}{%
    \noalign {\ifnum 0=`}\fi \hrule height 0.8pt
    \futurelet \reserved@a \@xhline
}
\newcolumntype{"}{@{\hskip\tabcolsep\vrule width 0.8pt\hskip\tabcolsep}}
\makeatother

\title{Restart FISTA with Global Linear Convergence}

\author{Teodoro Alamo, Pablo Krupa and Daniel Limon
\thanks{T. Alamo, P. Krupa and D. Limon are at the Systems Engineering and Automation Department, University of Seville, Spain. E-mail: {\tt \small \{talamo,pkrupa,dlm\}@us.es} }
\thanks{The authors acknowledge MINERCO and FEDER funds for funding project
	DPI2016-76493-C3-1-R, and MCIU and FSE for the FPI-2017 grant. }
\thanks{This paper constitutes an extended and revised version of \cite{AlamoECC:19}. Some of the technical results presented in this paper are used in \cite{AlamoCDC:19}.}
}

\begin{document}
\maketitle
\pagestyle{plain}

\begin{abstract}
Fast Iterative Shrinking-Threshold Algorithm (FISTA) is a popular fast gradient descent method (FGM) in the field of large scale convex optimization problems. However, it can  exhibit undesirable periodic oscillatory behaviour in some applications that slows its convergence. Restart schemes seek to improve the convergence of FGM algorithms by suppressing the oscillatory behaviour. Recently, a restart scheme for FGM has been proposed that provides linear convergence for non strongly convex optimization problems that satisfy a quadratic functional growth condition. However, the proposed algorithm requires prior knowledge of the optimal value of the objective function or of the quadratic functional growth parameter. In this paper we present a restart scheme for FISTA algorithm, with global linear convergence, for non strongly convex optimization problems that satisfy the quadratic growth condition without requiring the aforementioned values. We present some numerical simulations that suggest that
the proposed approach outperforms other restart FISTA schemes.  
\end{abstract}

\subsubsection*{Keywords}
Fast gradient method, restart FISTA, convex optimization, linear convergence, quadratic functional growth condition.

\section{Introduction} \label{sec:Intro}
Fast gradient methods (FGM) were introduced by Yurii Nesterov in \cite{Nesterov83}, \cite{Nesterov04}, where it was shown that these methods provide a convergence rate \textit{O}$(1/k^2)$ for smooth convex optimization problems with non strongly convex objective functions \cite{Nesterov04}, where $k$ is the iteration counter. These methods were generalized to composite non smooth convex optimization problems in \cite{Beck09}, \cite{Nesterov13}, \cite{Tseng:08}. The resulting algorithm is commonly known as FISTA algorithm \cite{Beck09}. Because of its complexity certification, it is often used in the context of embedded model predictive control \cite{koegel2011fast}, \cite{Richter12}, \cite{Krupa:18}. Another possibility to address composite convex optimization problems is to use splitting methods like ADMM \cite{boyd2011distributed}, \cite{Giselsson:TAC:17}, \cite{Banjac:TAC:18}.

FISTA algorithms can be applied in a primal setting (as in the Lasso problem \cite{Beck09}), or in a dual one \cite{Richter:13}, \cite{Beck2014Dual}. They can be thought of as a momentum method, since the linearization point at each iteration depends on the previous iterations. Since the momentum grows with the iteration counter, the algorithm can exhibit undesirable periodic oscillating behavior for certain applications, which slows the convergence rate. To mitigate this, restart schemes have been proposed in the literature which stop the algorithm when a certain criteria is met. It is then restarted using the last value provided by the stopped algorithm as the new initial condition \cite{Donoghue:13}, \cite{Alamir:13}, \cite{Giselsson:14CDC}.

In \cite{Donoghue:13} two heuristic restart schemes for FGM are proposed which exhibit improved convergence rates over non-restart FGM schemes. These restart schemes reset the momentum of the FGM in order to eliminate the undesirable oscillations whenever the periodical behavior is detected. A restart scheme similar to the ones in \cite{Donoghue:13} with \textit{O}$(1/k^2)$ convergence rate for smooth convex optimization is presented in \cite{Giselsson:14CDC}. In \cite{Kim:18}, an algorithm is proposed that uses the restart schemes from \cite{Donoghue:13}. Numerical results show improvements over previous restart schemes for FGM, but no theoretical results on convergence rates are provided.

Recently, linear convergence rate has been derived for several first order methods applied to convex optimization problems with non strongly convex objective functions that satisfy a relaxation of the strong convexity known as the quadratic functional growth \cite{Necoara:18}.

In \cite[Subsection 5.2.2]{Necoara:18} a restarting scheme of FGM is presented with global linear convergence rate for convex optimization problems that satisfy the functional growth condition with parameter $\mu$. However, in order to implement this strategy, prior knowledge is needed of either the optimal value of the objective function or the value of $\mu$, which can be challenging to compute.

In this paper we propose a novel restart scheme for FISTA algorithm applied to solving convex constrained problems.
We show that the algorithm guarantees global linear convergence rate \textit{O}$(1/\sqrt{\mu})$ for convex optimization problems with non strongly convex objective functions that satisfy the quadratic functional growth condition with parameter $\mu$. The proposed algorithm does not require prior knowledge of the value of $\mu$ or of the optimal value of the objective function. We provide theoretical upper bounds on the number of iterations of the algorithm needed to achieve a given accuracy.

Additionally, we show numerical results comparing the proposed algorithm with the heuristic restart schemes from \cite{Donoghue:13} and the restart scheme from \cite{Necoara:18} for Lasso problems. 

In Section \ref{sec:Formulation} we introduce the problem formulation. Section \ref{sec:FISTA} presents FISTA algorithm and some restart schemes. The convergence rate of non restart FISTA algorithm under the satisfaction of the quadratic functional growth condition is presented in Section \ref{sec:Conv_FISTA}. In Section \ref{sec:OurFISTA} we present the proposed restart scheme for FISTA and state its global linear convergence. Numerical results comparing the proposed algorithm with other restart schemes applied to FISTA are shown in Section \ref{sec:Results}. Finally, conclusions are presented in Section \ref{sec:Conclusions}.

\subsubsection*{Notation}
Given vectors $x$ and $y$, we denote by $\sp{x}{y}$ their scalar product, i.e. $\sp{x}{y}\doteq x\T y$. Given vector $x$, $\|x\|_2$ denotes its Euclidean norm ($\|x\|_2\doteq \sqrt{x\T x}$), and $\| \cdot \|_1$ denotes its $l_1$-norm (sum of the absolute values of the components of $x$). Given $R\succ 0$ we denote by $\| \cdot \|_R$ the weighted Euclidean norm $\| x \|_R \doteq \sqrt{x\T R x}$, and by $\| x \|_* \doteq \| x \|_{R^{-1}}$ its dual norm. $\ln(\cdot)$ is the natural logarithm and $e$ is Euler's number. $\lfloor x\rfloor$ denotes the largest integer smaller than or equal
to $x$; $\lceil x \rceil $ denotes the smallest integer greater than or equal to $x$. 
Given a set $\cX\subseteq \R^n$ we denote by $I_{\cX}$ its indicator function. That is, ${I_{\cX}(x)=0}$ if ${x\in \cX}$, and ${I_{\cX}(x)=\infty}$ if $x\not\in \cX$. The relative interior of set $\cX$ is denoted by $\ri{\cX}$. 
Given the extended real valued function ${f:\R^n\to (-\infty,\infty]}$ we denote by $\dom{f}$  its effective domain. That is, ${ \dom{f}\doteq\set{x\in \R^n}{f(x)<\infty}}$. We denote by $\epi{f}$ the epigraph of $f$. That is, ${\epi{f}\doteq\set{(x,t)\in \R^n\times \R}{f(x)\leq t}}$. We say that function ${f:\R^n\to (-\infty,\infty]}$ is closed if its epigraph is a closed set. 
We say that ${f:\R^n\to (-\infty,\infty]}$ is proper if its effective domain is not empty. That is, if $f$ is not identically equal to $\infty$. We say that a vector ${d\in \R^n}$ is a subgradient of $f$ at a point ${x\in\dom{f}}$ if ${f(y)\geq f(x)+\sp{d}{y-x}}$, ${\forall y\in \R^n}$. The set of all subgradients of $f$ at $x$ is called the subdifferential of $f$ at $x$ and is denoted by $\partial f(x)$.

\section{Problem Formulation} \label{sec:Formulation}

We address the problem of solving the composite convex minimization problem
\begin{equation} \label{eq:OP}
    f^* = \min\limits_{x\in \cX} f(x) = \min\limits_{x\in \cX} \Psi(x) + h(x),
\end{equation}
under the following assumption.

\begin{assumption}\label{assum:conv:smooth}
	We assume that
	\blista
	\item $h:\R^n\to \R$ is a smooth differentiable convex function. That is, there is $R\succ 0$ such that the inequality
	\begin{equation}\label{ineq:smooth:R}
	    h(x) \leq h(y)+\sp{\dhh(y)}{x-y}+\frac{1}{2}\|x-y\|_R^2,
	\end{equation}
	is satisfied for every $x\in \R^n$ and $y\in \R^n$.
	\item $\Psi:\R^n\to (-\infty,\infty]$ is a closed convex function and $\cX\subseteq \R^n$ is a closed convex set.
	\item Denote $f\doteq \Psi+h$. The minimization problem
	$$\min\limits_{x\in \cX} f(x)$$
	is solvable. That is, there is $x^* \in \cX \bigcap\dom{\Psi}$ such that $f^*=f(x^*)=\inf\limits_{x\in \cX} f(x)$.
	\elista
\end{assumption}

We notice that it is standard to write down the first point of Assumption \ref{assum:conv:smooth} as
\begin{equation}\label{equ:h:with:L}
    h(x) \leq h(y)+\sp{\dhh(y)}{x-y}+\frac{L}{2}\|x-y\|_{S}^2,
\end{equation}
where parameter $L$ serves to characterize the smoothness of $h$ and $S$ is a positive definite matrix. Constant $L$ provides a bound on the Lipschitz constant of the gradient $\nabla h(\cdot)$ \cite[Subsection 2.1]{Nesterov04}. Since $$\frac{L}{2}\|x-y\|_{S}^2 = \frac{1}{2} \|x-y\|_{  L S}^2,$$ we have that (\ref{equ:h:with:L}) implies (\ref{ineq:smooth:R}) if we take $R=LS$. This simplifies the algebraic expressions needed to analyze the convergence of the proposed algorithm.

We notice that Assumption \ref{assum:conv:smooth} guarantees that the minimization problem (\ref{eq:OP}) is solvable. The optimal set $\Omega$ is defined as
$$ \Omega \doteq \set{x}{x\in \cX, f(x)=f^*}.$$
This set is a singleton if $f(x)$ is strictly convex. Given ${x\in\R^n}$ we will denote $\bx$ its closest element in the optimal set $\Omega$ (with respect to the norm $\|\cdot\|_R$). That is,
\begin{equation}\label{def:bx}
    \bx \doteq  \arg \min\limits_{z\in\Omega} \|x-z\|_R.
\end{equation}
Given $y\in \R^n$, one could use the local information given by $\dhh(y)$ to minimize the value of $f=\Psi+h$ around $y$.
Under Assumption \ref{assum:conv:smooth}, this can be done obtaining the minimizer of the strictly convex optimization problem
$$ \min\limits_{x\in \cX} \; \Psi(x) +\sp{\dhh(y)}{x-y}+\frac{1}{2} \|x-y\|_R^2.$$
It is well known that this problem is solvable and has a unique solution if Assumption \ref{assum:conv:smooth} holds (see, for example, Subsection 6.1 in \cite{beck2017} for an analogous result). For completeness we provide a proof of this statement in Appendix \ref{appen:existence:uniqueness} (see Property \ref{Prop:solvable:unique:grad}). 

The solution to this optimization problem leads to the notion of composite gradient mapping \cite{Nesterov13}, which constitutes a generalization of the gradient mapping that can be found in \cite[Subsection 2.2]{Nesterov04} for the particular case $\Psi(\cdot)=0$. See also \cite{Beck09} for the particular case $\cX=\R^n$.
\begin{definition}[Composite Gradient Mapping $g(y)$]\label{def:composite:gradient:mapping}
~\\Under Assumption \ref{assum:conv:smooth}, and given $y\in \R^n$, we define
\begin{align*}
    y^+ &\doteq \arg \min\limits_{x\in \cX} \; \Psi(x)+ \sp{\dhh(y)}{x-y}+\frac{1}{2} \|x-y\|_R^2,\\
    g(y) &\doteq R(y-y^+).
\end{align*}
\end{definition}
We notice that the composite gradient mapping is closely related to the notion of proximal operator \cite{Parikh13}, \cite[Chapter 6]{beck2017}. For example, one could state, after some manipulations, the computation of the composite gradient mapping as the computation of a proximal operator. In the context of optimal gradient methods, it is assumed that the computation of $y^+$ is cheap. This is the case when $\cX$ is a simple set (box, $\R^n$, etc.), $R$ diagonal, and $\Psi(\cdot)$ a separable function. For example, in the well known Lasso optimization problem, the computation of $y^+$ resorts to the computation of the shrinkage operator \cite{Beck09}. See \cite{Combettes11}, Section 6 of \cite{Parikh13}, Chapter 28 in \cite{Bauschke:11}, or Chapter 6 in \cite{beck2017}, for numerous examples in which the computation of the composite gradient mapping is simple.

The following property gathers well-known properties of the composite gradient mapping $g(y)$ and its dual norm $\| g(y) \|_*=\|g(y)\|_{R^{-1}}$ \cite{Beck09}, \cite{Nesterov13}. For completeness, we
include the proof in Appendix \ref{appen:proj:gradient}.

\begin{property}\label{prop:proj:gradient}
Suppose that Assumption \ref{assum:conv:smooth} holds.
Then,
\blista
\item For every $y\in \R^n$ and $x\in \cX$:
\begin{align*}
 f(y^+) -f(x) &\leq \sp{g(y)}{y^+-x}  + \frac{1}{2}\|g(y)\|_*^2\\
              & = \sp{g(y)}{y-x}  - \frac{1}{2}\|g(y)\|_*^2 \\
              & = - \frac{1}{2}\| y^+ -x \|^2_R + \frac{1}{2}\|y-x\|_R^2.
\end{align*} \label{prop:proj:gradient_1}
\item For every $y\in \cX$: $$ \frac{1}{2} \| g(y)\|_*^2 \leq f(y)-f(y^+) \leq f(y)-f^*.$$ \label{prop:proj:gradient_2}
\elista
 \end{property}

The composite gradient serves to characterize optimality \cite{Nesterov13}. That is, under Assumption \ref{assum:conv:smooth} we have the following equivalence
$$ y\in \Omega \Leftrightarrow  g(y)=0. $$
This fact is proved in Appendix \ref{appen:optimaliticharacterization}. 
   
 \section{Restart FISTA Schemes } \label{sec:FISTA}

For a given initial condition $z\in \R^n$, a minimum number of iterations $k_{min} \geq 0$, and an exit condition $E_c$, the non restart FISTA algorithm \cite{Beck09} is shown in Algorithm \ref{alg:FISTA}. This algorithm solves 
$\min\limits_{x\in\cX} \; h(x)+\Psi(x)$ under Assumption \ref{assum:conv:smooth}. 

\begin{algorithm}
	\DontPrintSemicolon
	\caption{FISTA} \label{alg:FISTA}
	\Require{$z \in \R^n$, $k_{min} \geq 0$, $E_c$}
	$y_0 = x_0=z^+$, $t_0 = 1$, $k = 0$\;
	\Repeat{$E_c$ and $k \geq k_{min}$} {
		$k = k + 1$\;
        $x_{k} = y_{k-1}^+$ \label{alg:step:grad}\;
		$t_k = \fracg{1}{2} \left( 1 + \sqrt{1 + 4 t^{2}_{k-1}}\, \right)$ \label{FISTA:tk}\;
		$y_k = x_k + \fracg{t_{k-1} - 1}{t_k} (x_k - x_{k-1})$\;
	    Compute exit condition $E_c$\;}
	\KwOut{$r = x_k$, $n= k$}
\end{algorithm}

Since the optimality of $x_k$ is equivalent to $g(x_k)=0$ (see Property \ref{prop:charac:optimality} in Appendix \ref{appen:optimaliticharacterization}), a typical choice for non restart FISTA schemes is to choose $k_{min}$ equal to zero and codify the exit condition
\begin{equation*} \label{eq:exit} 
     \| g(x_k) \|_* \leq \epsilon,
\end{equation*}
where $\epsilon >0$ is an accuracy parameter. It is also common to use the exit condition $\|g(y_{k-1})\|_*\leq \epsilon$, since this exit condition requires $y_{k-1}^+$, which has already been computed in step \ref{alg:step:grad} of the algorithm.

It is well known that under Assumption \ref{assum:conv:smooth}, see also (\ref{equ:h:with:L}), the iterations of non-restart FISTA satisfy \cite{Beck09,Nesterov13}, 
\begin{equation} 
    f(x_k)-f^* \leq \frac{2}{(k+1)^2}\| x_0-\bx_0\|^2_R, \; \forall k\geq 1, \label{equ:conv:FISTA}
\end{equation}
where $\bar{x}_0$ represents the point in the optimal set $\Omega$ closest to the initial condition $x_0$ of the algorithm (see (\ref{def:bx})). For the sake of completeness, we present a detailed proof of this claim in Appendix \ref{app:conver:Fista}. 
We also prove in the same appendix that the sequence $\{y_k\}$ generated by Algorithm \ref{alg:FISTA} (FISTA) satisfies
$$ \|g(y_k)\|_* \leq \fracg{4\| x_0-\bx_0\|_R}{k+2}, \forall k\geq 0.$$  

In restart schemes, one invokes several times FISTA algorithm with a relaxed exit condition. Typical choices are (see \cite{Donoghue:13}),
\blista
    \item Function scheme: 
    \begin{equation} \label{eq:Functional}
        E_c^{f} = \text{True} \Leftrightarrow f(x_{k})\geq f(x_{k-1}).
    \end{equation}
    \item  Gradient scheme: 
    \begin{equation} \label{eq:Gradient}
        E_c^g = \text{True} \Leftrightarrow \langle g(y_{k-1}), x_{k-1} - x_{k} \rangle \leq 0.
    \end{equation}
    \elista

Given initial condition $r_0\in \cX$, a minimum number of iterations $k_{min} \geq 0$, an exit condition $E_c$, and an accuracy parameter $\epsilon > 0$, the standard restart FISTA algorithm is shown in Algorithm \ref{alg:RestartFISTA}.

\begin{algorithm}
    \DontPrintSemicolon
    \caption{Restart FISTA} \label{alg:RestartFISTA}
    \Require{$r_0 \in \mathcal X$, $k_{min} \geq 0$, $\epsilon >0$, $E_c$}
    $j = 0$\;
    \Repeat{$\vert\vert g(r_j)\vert\vert_* \leq \epsilon$}{
        $j = j + 1$\;
        $r_{j} = \text{FISTA}(r_{j-1}, k_{min}, E_c)$\;
    } 
    \KwOut{$x^* = r_j$}
\end{algorithm}

The implementation of Algorithm \ref{alg:RestartFISTA} usually provides better performance than the original non restart version \cite{Donoghue:13}, \cite{Giselsson:14CDC}.

\section{Convergence of Restart FISTA under a quadratic functional growth condition} \label{sec:Conv_FISTA} 

It has been recently shown in \cite{Necoara:18} that some relaxations of the strong convexity conditions of the objective function are sufficient for obtaining linear convergence for several first order methods. In particular, the following relaxation of strong convexity suffices to guarantee linear convergence of different gradient optimization schemes for smooth functions ($\Psi(\cdot)=0$). See \cite[Subsection 5.2.2]{Necoara:18}. 

\begin{assumption}[Quadratic Functional Growth]\label{assump:quadratic:growth}
We assume that the optimization problem 
$$ f^* = \min\limits_{x\in \cX}\, f(x) $$  
is solvable and satisfies the following quadratic functional growth condition with parameter $\mu>0$:
$$ f(x) - f^* \geq \frac{\mu}{2}\| x-\bx\|_R^2, \; \forall x\in \cX, $$
where $\bx$ denotes the closest element to $x$ in the optimal set $\Omega$ (see (\ref{def:bx})).
\end{assumption}

As can be seen in \cite[Subsection 3.4]{Necoara:18}, strong convexity implies quadratic functional growth. This means that the quadratic functional growth setting encompasses a broad family of convex functions. 

It is also shown in \cite[Subsection 5.2.2]{Necoara:18} that if the value of $f^*$ is known and $\Psi(\cdot)=0$, then a restart FISTA based on the exit condition 
\begin{equation} \label{E:Optimal}
    E_c^* = \text{True} \Leftrightarrow f(x_{k}) - f^* \leq \frac{f(x_0) - f^*}{e^2},
\end{equation}
exhibits global linear convergence. This exit condition is easily implementable if the optimal value $f^*$ is known. This is the case, for example, in some formulations of feasibility optimization problems, in which the optimal value $f^*$ is equal to zero for every feasible solution. This restart scheme corresponds to an optimal restart rate of $\frac{2e}{\sqrt{\mu}}$  \cite[Subsection 5.2.2]{Necoara:18}. 

We present now a novel result that further characterizes the convergence properties of the non restart FISTA algorithm under Assumption \ref{assump:quadratic:growth}. 

\begin{property}\label{prop:FISTA:EC}
Under Assumptions \ref{assum:conv:smooth} and \ref{assump:quadratic:growth}, the iterations of FISTA algorithm satisfy
\blista
\item $f(x_k)-f^* \leq \fracg{4 (f(x_0)-f^*)}{\mu(k+1)^2}$, for all $k\geq 1$. \label{prop:FISTA:EC_1}
\item $f(x_k) \leq f(x_0)$, for all $k\geq \left\lfloor \frac{2}{\sqrt{\mu}}\right\rfloor$. \label{prop:FISTA:EC_2}
\item $f(x_k)-f^* \leq  \fracg{f(x_0)-f(x_k)}{e}$, for all $k\geq \left\lfloor \frac{2\sqrt{e+1}}{\sqrt{\mu}}\right\rfloor$. \label{prop:FISTA:EC_3}
\elista
\end{property}
\vspace{.5em}
\proof See Appendix \ref{appen:FISTA:EC}. 

\section{Restart FISTA with global linear convergence} \label{sec:OurFISTA}

In this section we propose a novel restart FISTA algorithm (Algorithm \ref{alg:OurRestartFISTA}) that exhibits global linear convergence under the quadratic functional growth condition. 
The algorithm uses exit condition $E_c^l$, which is defined to be true if the following two conditions are satisfied,
\begin{subequations}
\label{E:Our}
    \begin{align}[left={E_c^l = \text{True} \Leftrightarrow  \empheqlbrace}]
    f(x_m) - f(x_k) & \leq \frac{f(x_0) - f(x_m)}{e} \label{E:Our_e}\\
    f(x_k) & \leq f(x_0), \label{E:Our_f}
  \end{align}
\end{subequations}
with $m = \lfloor\frac{k}{2}\rfloor +1$.

\begin{algorithm}
	\DontPrintSemicolon
    \caption{Linearly Convergent Restart FISTA (LCR-FISTA)} \label{alg:OurRestartFISTA}
    \Require{$r_0 \in \mathcal X$, $\epsilon >0$}
    $n_0 = 0$, $j = 1$\;
    $[r_1, n_1] = \text{FISTA}(r_0, n_0, E_c^l)$\;
    \Repeat{$\vert\vert g(r_j)\vert\vert_* \leq \epsilon$}{
        $j = j + 1$\;
        $[r_{j}, n_{j}] = \text{FISTA}(r_{j-1}, n_{j-1}, E_c^l)$\;
        \If{$f(r_{j-1}) - f(r_{j}) > \fracg{1}{e}\left(f(r_{j-2}) - f(r_{j-1})\right)$}{
            $n_j=2n_{j-1}$\; 
        }
    } 
    \KwOut{$r^* = r_j$}
\end{algorithm}

Inequality (\ref{E:Our_f}) guarantees that the output of the FISTA algorithm is no larger than the one corresponding to its initial condition.

As it is stated in the following property, one of the main features of the proposed algorithm is that the number of iterations $n_j$ required at each FISTA iteration ${[r_j,n_j]=FISTA(r_{j-1},n_{j-1},E_c^l)}$ is upper bounded by ${\frac{4\sqrt{e+1}}{\sqrt{\mu}} \approx \frac{7.72}{\sqrt{\mu}}}$. Moreover, the number of iterations required by the proposed algorithm to attain a given accuracy $\epsilon$ is upper bounded by $$\frac{16}{\sqrt{\mu}} \left\lceil \ln\left(1+\frac{2(f(r_0)-f^*)}{\epsilon^2}\right)\right\rceil.$$

\vspace{5em}

\begin{property}\label{prop:L:Convergence:FISTA} 
Suppose that Assumptions \ref{assum:conv:smooth} and \ref{assump:quadratic:growth} hold. Then, the sequences $\{r_j\}$, $\{n_j\}$  provided by Algorithm \ref{alg:OurRestartFISTA} satisfy
\blista
\item $\fracg{1}{2}\| g(r_{j-1})\|_*^2 \leq f(r_{j-1})-f(r_j)$, $\forall j\geq 1$.
\vspace{.3em}
\item $ n_j \leq \fracg{4\sqrt{e+1}}{\sqrt{\mu}}$, $\forall j\geq 0$. \label{prop:L:Convergence:FISTA_1}
\vspace{.3em}
\item The number of iterations ($\Sum{i=0}{j} n_i$) required to guarantee $\| g(r_j)\|_*\leq \epsilon$ is no larger than 
$$ \fracg{16}{\sqrt{\mu}}\left\lceil \ln\left(1+\frac{2(f(r_0)-f^*)}{\epsilon^2}\right)\right\rceil.$$ \label{prop:L:Convergence:FISTA_4} 
\elista
\end{property}

\proof See Appendix \ref{appen:L:Convergence:FISTA}.

We notice that the factor 16 in the worst case complexity analysis is conservative. The authors claim that a better factor might be obtained at the expense of a more involved proof. 

\section{Numerical results} \label{sec:Results}

We consider a weighted Lasso problem of the form
\begin{equation} \label{eq:Lasso}
    \min\limits_{x} \frac{1}{2 N} \| A x - b \|^2_2 + \| W x \|_1,
\end{equation}
where $x \inR{n}$, $A \inR{N \times n}$ is sparse with an average of $90\%$ of its entries being zero (sparsity was generated by setting a $0.9$ probability for each element of the matrix to be $0$), $n > N$, and $b \inR{N}$. Each nonzero element in $A$ and $b$ is obtained from a Gaussian distribution with zero mean and covariance 1. $W \inR{n \times n}$ is a diagonal matrix with elements obtained from a uniform distribution on the interval $[0, \alpha]$. 

We note that Lasso problems (\ref{eq:Lasso}) can be reformulated in such a way that they satisfy the quadratic growth condition \cite[Section 6.3]{Necoara:18}. For this problem, inequality (\ref{ineq:smooth:R}) of Assumption \ref{assum:conv:smooth} is satisfied, for instance, for a matrix $R$ chosen as
\begin{equation*}
    R_{i,i} = \sum_{j=1}^{n} \vert H_{i,j} \vert,
\end{equation*}
with $H = \frac{1}{N} A\T A$. This is due to the Gershgorin Circle Theorem \cite[Subsection 7.2]{Golub96}. See also \cite[Section 6]{Nesterov13}.

We show the results of applying algorithms \ref{alg:RestartFISTA} and \ref{alg:OurRestartFISTA} with an accuracy parameter $\epsilon = 10^{-11}$ using different restart schemes and values of $N$, $n$ and $\alpha$. We take $r_0 = 0$.

The restart schemes shown are $E_c^f$ (\ref{eq:Functional}) and $E_c^g$ (\ref{eq:Gradient}) from \cite{Donoghue:13}, restart condition $E_c^*$ (\ref{E:Optimal}) \cite{Necoara:18}, and the restart condition $E_c^l$ (\ref{E:Our}) proposed in this paper (using Algorithm \ref{alg:OurRestartFISTA}). Additionally, we show the results of applying FISTA algorithm without using a restart scheme. In order to provide a fair comparison between the performance of the restart schemes, the algorithms are exited as soon as a value of $y_k$ that satisfies $\| g(y_{k-1}) \|_* \leq \epsilon$ is found. We note that, in order to implement the restart scheme based on $E_c^*$, we had to previously compute the optimal value $f^*$, which was done by using Algorithm \ref{alg:OurRestartFISTA} with $\epsilon = 10^{-12}$.

Tables \ref{tab:Test1} to \ref{tab:Test3} show results of performing $100$ tests with different randomized problems (\ref{eq:Lasso}) that share common values of parameters $N$, $n$ and $\alpha$. Tables show the average, median, maximum and minimum number of iterations.

{\renewcommand{\arraystretch}{1.4}%
\begin{table}[ht]
    \centering
    \caption{Test 1. Comparison between restart schemes} 
    \label{tab:Test1}
    \begin{threeparttable}
        \begin{tabular}{|c"c|c|c|c|c|}
        \hline
        Exit Cond.         & $E_c^l$   & No restart & $E_c^f$   & $E_c^g$   & $E_c^*$   \\\thickhline
        Avg. Iter.         & $670.6$   & $8207.2$   & $1648.7$  & $687.5$   & $1569.5$  \\\hline
        Median Iter.       & $676$     & $8241$     & $1608.5$  & $666.5$   & $1571$    \\\hline
        Max. Iter.         & $783$     & $10109$    & $2156$    & $930$     & $2053$    \\\hline
        Min. Iter.         & $570$     & $6737$     & $1192$    & $567$     & $917$     \\\hline
        \end{tabular}
        \begin{tablenotes}[flushleft] \footnotesize
            \item Results of $100$ tests with $N=600$, $n=800$, $\alpha = 0.01$, $\epsilon = 10^{-11}$.
        \end{tablenotes}
    \end{threeparttable}
\end{table}}

{\renewcommand{\arraystretch}{1.4}%
\begin{table}[ht]
    \centering
    \caption{Test 2. Comparison between restart schemes} 
    \label{tab:Test2}
    \begin{threeparttable}
        \begin{tabular}{|c"c|c|c|c|c|}
        \hline
        Exit Cond.         & $E_c^l$  & No restart & $E_c^f$   & $E_c^g$   & $E_c^*$         \\\thickhline
        Avg. Iter.         & $1683.7$ & $34116.4$  & $7743.3$  & $1606.7$  & $4601.9$        \\\hline
        Median Iter.       & $1659$   & $33127.5$  & $7242$    & $1594$    & $4503$          \\\hline
        Max. Iter.         & $2162$   & $51201$    & $14080$   & $2201$    & $7266$          \\\hline
        Min. Iter.         & $1406$   & $24539$    & $3894$    & $1306$    & $2499$          \\\hline
        \end{tabular}
        \begin{tablenotes}[flushleft] \footnotesize
            \item Results for $100$ tests with $N=600$, $n=800$, $\alpha = 0.003$, $\epsilon = 10^{-11}$.
        \end{tablenotes}
    \end{threeparttable}
\end{table}}

{\renewcommand{\arraystretch}{1.4}%
\begin{table}[ht]
    \centering
    \caption{Test 3. Comparison between restart schemes} 
    \label{tab:Test3}
    \begin{threeparttable}
        \begin{tabular}{|c"c|c|c|c|c|}
        \hline
        Exit Cond.         & $E_c^l$   & No restart & $E_c^f$   & $E_c^g$   & $E_c^*$   \\\thickhline
        Avg. Iter.         & $705.9$   & $8379.5$   & $1786.3$  & $686$     & $1709.4$  \\\hline
        Median Iter.       & $704.5$   & $8135.5$   & $1773$    & $680.5$   & $1703$    \\\hline
        Max. Iter.         & $873$     & $12055$    & $3218$    & $892$     & $2512$    \\\hline
        Min. Iter.         & $547$     & $5943$     & $987$     & $529$     & $1042$    \\\hline
        \end{tabular}
        \begin{tablenotes}[flushleft] \footnotesize
            \item Results for $100$ tests with $N=300$, $n=400$, $\alpha = 0.01$, $\epsilon = 10^{-11}$.
        \end{tablenotes}
    \end{threeparttable}
\end{table}}

\begin{figure}[ht]
    \centering
    \includegraphics[width=\columnwidth]{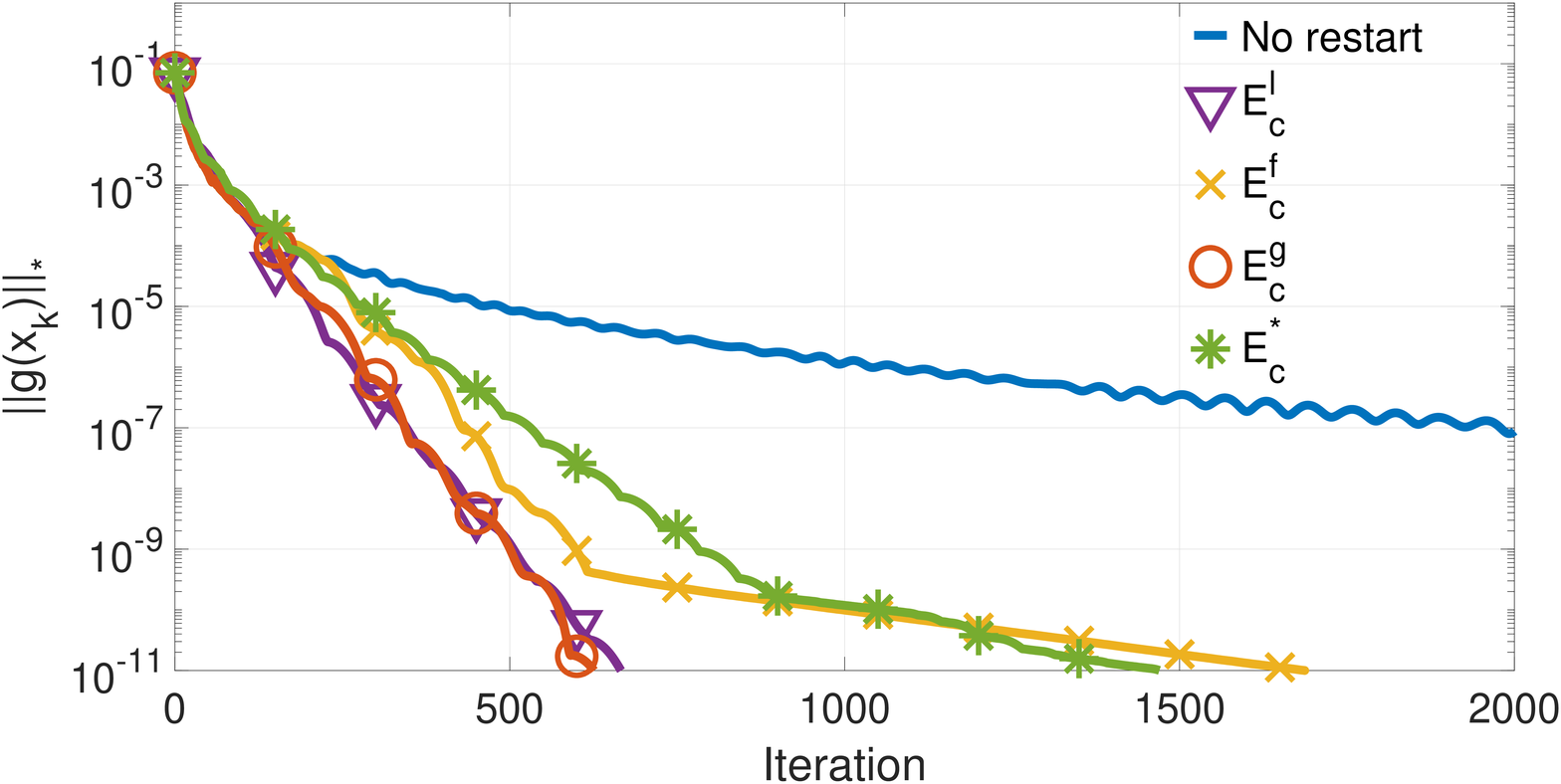}
    \caption{Value of $\| g(y_k) \|_*$ for a problem (\ref{eq:Lasso}) of Test 1.}
    \label{fig:Test1_gx}
\end{figure}

\begin{figure}[ht]
    \centering
    \includegraphics[width=\columnwidth]{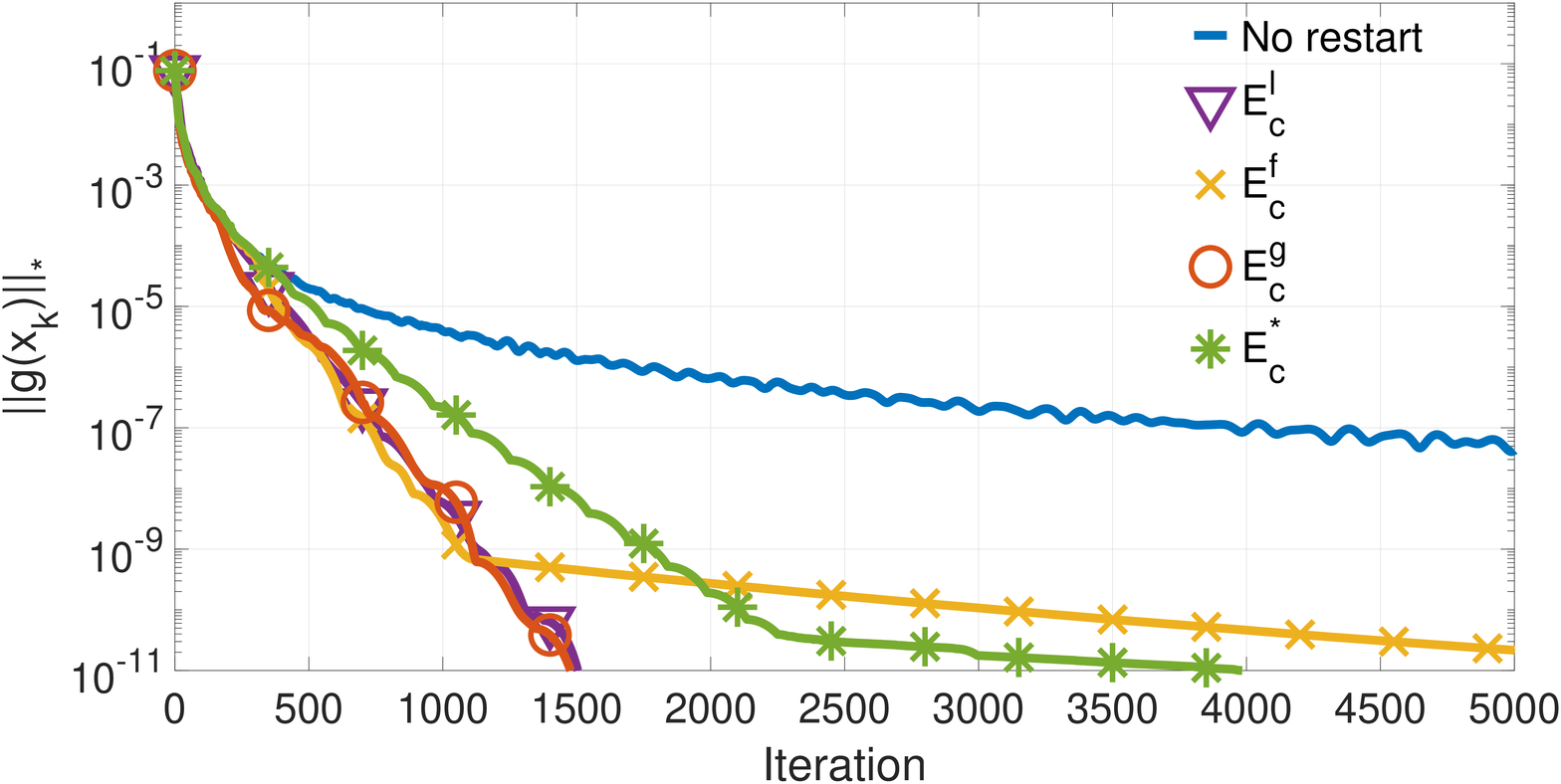}
    \caption{Value of $\| g(y_k) \|_*$ for a problem (\ref{eq:Lasso}) of Test 2.}
    \label{fig:Test2_gx}
\end{figure}

\begin{figure}[ht]
    \centering
    \includegraphics[width=\columnwidth]{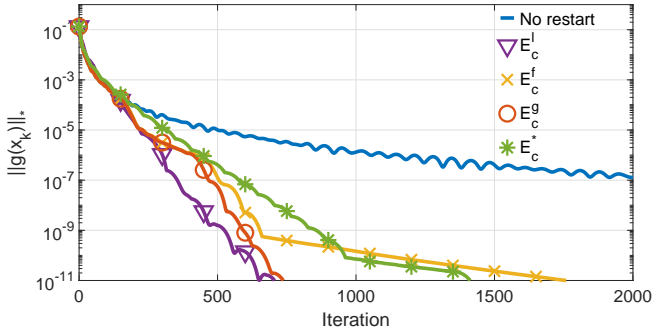}
    \caption{Value of $\| g(y_k) \|_*$ for a problem (\ref{eq:Lasso}) of Test 3.}
    \label{fig:Test3_gx}
\end{figure}

Figures \ref{fig:Test1_gx} to \ref{fig:Test3_gx} show the value of $\| g(x_k) \|_*$ for a randomly selected problem out of the randomized problems used to compute the results shown in tables \ref{tab:Test1} to \ref{tab:Test3}, respectively.

\begin{figure}[ht]
    \centering
    \includegraphics[width=\columnwidth]{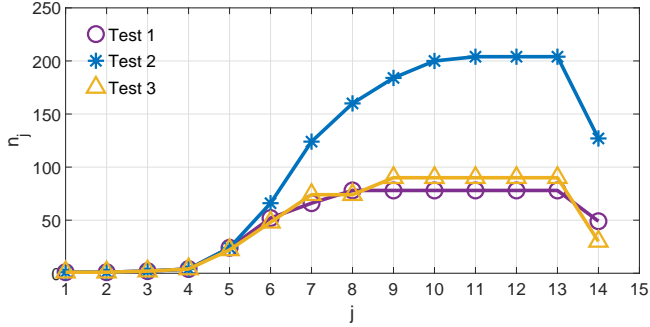}
   \caption{Value of $n_j$ obtained using Algorithm \ref{alg:OurRestartFISTA} for the problems (\ref{eq:Lasso}) whose result are shown in Figures \ref{fig:Test1_gx} to \ref{fig:Test3_gx}.}
    \label{fig:nj}
\end{figure}

Figure \ref{fig:nj} shows the value of $n_j$ at each iteration $j$ of Algorithm \ref{alg:OurRestartFISTA} for the three examples whose results are shown in Figures \ref{fig:Test1_gx} to \ref{fig:Test3_gx}. Note that the final value of $n_j$ is lower than the previous one in all three instances due to the algorithm exiting as soon as the condition $\| g(y_{k-1}) \|_* \leq \epsilon$ is satisfied.

\section{Conclusions} \label{sec:Conclusions}

In this paper we have presented a novel restart scheme with guaranteed global linear convergence. The algorithm relies on a quadratic functional growth condition. One of the advantages of the proposed algorithm is that it does not require the knowledge of the parameter $\mu$ that characterizes the quadratic functional growth condition, or the optimal value of the minimization problem. We provide an upper bound of the required number of iterations equal to $$\frac{16}{\sqrt{\mu}} \left\lceil \ln\left(1+\frac{2(f(r_0)-f^*)}{\epsilon^2}\right)\right\rceil.$$
 We have presented numerical evidence of the good performance of the algorithm when compared with other restarts schemes. It outperforms the restart scheme  based on the knowledge of the optimal value $f^*$.

\bibliography{Bib_Arxiv_Restart_FISTA_ECC19}


\begin{appendix}

\subsection{Existence and Uniqueness of Composite Gradient}\label{appen:existence:uniqueness}

We present in this appendix some well known facts about convex analysis that are required to analyze the properties of the composite gradient.

\begin{property}\label{prop:characterization:Psi:cX}
	Suppose that
	\blista
	\item $\Psi: \R^n \to (-\infty,\infty]$ is a closed convex function.
	\item $\cX \subseteq \R^n$ is a closed convex set.
	\item The set $\dom{\Psi}\bigcap\cX $ is non empty.
	\item $I_{\cX}:\R^n \to \{0,\infty\}$ is the indicator function of $\cX$. That is,
	$$ I_{\cX}(x) = \bsis{rl} 0 & \mbox{if } x\in \cX \\ \infty & \mbox{otherwise}. \esis$$
	\item The function $\Psi_{\cX}:\R^n\to (-\infty,\infty]$ is defined as
	$$\Psi_{\cX} (x)\doteq \Psi(x)+ I_{\cX}(x),\; \forall x\in \R^n.$$
	\elista
	Then
	\blista
	\item The function $\Psi_{\cX}$ is proper, closed, and convex.
    \item The relative interior of $\dom{\Psi_{\cX}}$ is non empty.
	\item There is $z\in \cX$ and $d\in \R^n$ such that $\Psi_{\cX}(z)<\infty$ and
	$$ \Psi_{\cX} (x) \geq \Psi_{\cX}(z) + \sp{d}{x-z},\;\; \forall x\in \R^n.$$
	\elista
\end{property}

\proof

From $\dom{\Psi}\bigcap\cX\neq \emptyset $ we have that both $\dom{\Psi}$ and $\cX$ are non empty.
The epigraph of the indicator function $I_{\cX}$ is, by definition,
\begin{eqnarray*}
    \epi{I_{\cX}} &=& \set{(x,t) \in \R^n \times \R }{I_{\cX}(x) \leq t }\\
    & =&  \set{(x,t) \in \R^n\times \R}{x\in \cX, 0 \leq t}.\end{eqnarray*}

Since $\cX$ and ${{\cal{T}}\doteq\set{t\in \R}{t\geq 0}}$ are non empty closed sets, ${\epi{I_{\cX}}=\cX\times {\cal{T}}}$ is also a non empty closed convex set. Thus, by definition,  ${I_{\cX}:\R^n\to \{0,\infty\}}$ is a closed convex function. Since both $\Psi$ and $I_{\cX}$ are closed convex functions, ${\Psi_{\cX}\doteq\Psi+I_{\cX}}$ is also a closed convex function (the sum of closed convex functions provides closed convex functions \cite[Proposition 1.1.5]{Bertsekas:09}). Since ${\dom{\Psi_{\cX}}=\dom{\Psi}\bigcap\cX\neq \emptyset}$, we infer that the domain of $\Psi_{\cX}$ is non empty. This implies that $\Psi_{\cX}$ is not identically equal to $\infty$. Moreover, since ${\Psi:\R^n\to (-\infty,\infty]}$ we have that ${\Psi_{\cX} : \R^n \to (-\infty,\infty]}$. We conclude that ${\Psi_{\cX}(x)>-\infty}$ for every ${x\in \R^n}$. From this and the fact that $\Psi_{\cX}$ is not identically equal to $\infty$ we have that $\Psi_{\cX}$ is proper.

Since $\dom{\Psi_{\cX}}$ is a non empty convex set, it has a non empty relative interior $\ri{\dom{\Psi_{\cX}}}$ (see \cite[Proposition 1.3.2]{Bertsekas:09}). 

It is a well know fact from convex analysis that the subdifferential of a proper convex function at a point in the relative interior of its domain is non empty \cite[Proposition 5.4.1]{Bertsekas:09}.  Suppose now that $z\in \ri{\dom{\Psi_{\cX}}}$. Since $\Psi_{\cX}$ is a proper convex function we have that the subdifferential of $\Psi_{\cX}$ at $z$ is non empty. This means, by definition, that there is $d\in \R^n$ such that
$$ \Psi_{\cX} (x) \geq \Psi_{\cX}(z) + \sp{d}{x-z},\;\; \forall x\in \R^n.$$
\QED

\begin{property}\label{Prop:solvable:unique:grad}
	Suppose that Assumption \ref{assum:conv:smooth} holds. Given any $y\in \R^n$, consider the quadratic function $h_y:\R^n \to \R$ defined as
	$$ h_y(x) \doteq \sp{\dhh(y)}{x-y}+\frac{1}{2} \|x-y\|^2_R.$$
	Then, the minimization problem
	\begin{equation}\label{opti:comp:grad:X}
	    \min\limits_{x\in \cX} \Psi(x)+h_y(x)
	\end{equation}
	is solvable and has a unique solution. That is, there exists a unique point $y^+\in \cX$ such that
	$$ \Psi(y^+)+h_y(y^+) =  \inf\limits_{x\in \cX} \Psi(x) + h_y(x) < \infty.$$
\end{property}

\proof

Notice that the minimization problem (\ref{opti:comp:grad:X}) is equivalent to
\begin{equation*}\label{opti:comp:grad:Indicator}
    \min\limits_{x\in \R^n} \Psi(x) + I_{\cX}(x) + h_y(x),
\end{equation*}
where $I_{\cX}$ is the indicator function of $\cX$. If we define ${\Psi_{\cX}\doteq \Psi+I_{\cX}}$ we can rewrite the original problem (\ref{opti:comp:grad:X}) as
\begin{equation*}\label{opti:comp:grad:Psi:cX}
    \min\limits_{x\in \R^n} \Psi_{\cX}(x)+ h_y(x).
\end{equation*}
We notice that the assumptions of Property \ref{prop:characterization:Psi:cX} are satisfied if Assumption $\ref{assum:conv:smooth}$ holds.
Thus, we infer from Property \ref{prop:characterization:Psi:cX} that $\Psi_{\cX}:\R^n\to (-\infty,\infty]$ is a proper closed convex function.
We also have that the quadratic function $h_y:\R^n \to \R$ is also proper and closed because it is a real valued continuous function (see \cite[Proposition 1.1.3]{Bertsekas:09}).
Since the sum of closed functions is closed (see \cite[Proposition 1.1.5]{Bertsekas:09}), we infer that $F_y\doteq\Psi_{\cX}+h_y$ is a closed function.
Moreover, from Property \ref{prop:characterization:Psi:cX} we also have that there is $z\in \cX$ and $d\in \R^n$ such that
\blista
    \item $\Psi_{\cX}(z)<\infty$.
    \item $ \Psi_{\cX} (x) \geq \Psi_{\cX}(z) + \sp{d}{x-z}$, $\forall x\in \R^n$.
\elista
Therefore,
\begin{equation}
    \begin{aligned} 
        F_y(z)&= \Psi_{\cX}(z)+h_y(z) = \gamma_z <  \infty, \label{equ:ineq:Psi:cX:hy}\\
        F_y(x)&=\Psi_{\cX}(x) + h_y(x) \\
        &\geq\Psi_{\cX}(z) + \sp{d}{x-z} + h_y(x),\; \forall x\in \R^n.
    \end{aligned}
\end{equation}
We infer from (\ref{equ:ineq:Psi:cX:hy}) that the closed function ${F_y:\R^n\to (-\infty,\infty]}$ is not identically equal to $\infty$ and therefore, proper.
We conclude that $F_y$ is a proper closed convex function. From Weiertrasss' Theorem (see Proposition 3.2.1 in \cite{Bertsekas:09}) we have that the set of minima of $F_y$ over $\R^n$ is nonempty and compact if there is a scalar $\bar{\gamma}$ such that the level set ${\Phi(\bar{\gamma})=\set{x}{F_y(x)\leq \bar{\gamma}}}$ is nonempty and bounded.
From (\ref{equ:ineq:Psi:cX:hy}) we have that $\Phi(\gamma_z)$ is nonempty. Moreover, we also infer from (\ref{equ:ineq:Psi:cX:hy}) that $\Phi(\gamma_z)$ is a bounded set because $F_y$ is lower bounded by a strictly convex quadratic function of $x$. We conclude that
\begin{eqnarray*}
	 \min\limits_{x\in \cX} \Psi(x)+h_y(x) &=& \min\limits_{x\in \R^n} \Psi_{\cX}(x)+ h_y(x) \\
	 &=& \min\limits_{x\in \R^n} F_y(x) \leq \gamma_z <\infty .
\end{eqnarray*}
is a solvable optimization problem. That is, there is $y^+\in \cX$ such that
$$ \Psi(y^+)+h_y(y^+) =  \inf\limits_{x\in \cX} \Psi(x) + h_y(x)  < \infty.$$
The set of minimizers consists of a single element $y^+$ because of the strictly convex nature of $F_y$ ($h_y$ is a strictly convex function).
\QED

\subsection{Proof of Property \ref{prop:proj:gradient}.}\label{appen:proj:gradient}

We prove in this appendix Property \ref{prop:proj:gradient}, which is rewritten here for the reader's convenience.

\begin{property}\label{prop:proj:in:appendix}
	Suppose that Assumption \ref{assum:conv:smooth} holds.
	Then,
	\blista
	\item For every $y\in \R^n$ and $x\in \cX$:
	\begin{subequations}
	   \begin{flalign}
	   f(y^+) -f(x) & \leq \sp{g(y)}{y^+-x}  + \frac{1}{2}\|g(y)\|_*^2 \label{ineq:cg:a}\\
	   &  = \sp{g(y)}{y-x}  - \frac{1}{2}\|g(y)\|_*^2 \label{ineq:cg:b}\\
	   & =  - \frac{1}{2}\| y^+ -x \|^2_R + \frac{1}{2}\|y-x\|_R^2. \label{ineq:cg:c}
	   \end{flalign}
	\end{subequations}
	\item For every $y\in \cX$: $$ \frac{1}{2} \| g(y)\|_*^2 \leq f(y)-f(y^+)\leq f(y)-f^*.$$
	\elista
\end{property}

\proof

From Property \ref{Prop:solvable:unique:grad} we have that there is a (unique) $y^+\in \cX$ such that
\begin{equation}\label{ineq:Psi:hy}
    \Psi(y^+) + h_y(y^+) \leq \Psi(x)+h_y(x), \; \forall x\in \cX,\end{equation}
where $ h_y(x) \doteq  \sp{\dhh(y)}{x-y}+\frac{1}{2} \|x-y\|^2_R$. 
Denote now $\Psi_{\cX}=\Psi+I_{\cX}$, where $I_{\cX}:\R^n\to \{0,\infty\}$ is the indicator function of $\cX$. Since $y^+\in \cX$ we have $I_{\cX}(y^+)=0$. Therefore, inequality
(\ref{ineq:Psi:hy}) implies
$$  \Psi_{\cX}(y^+) + h_y(y^+) \leq \Psi_{\cX}(x)+h_y(x), \; \forall x \in \R^n. $$
Denote now $F_y=\Psi_{\cX}+h_y$. From last inequality we have
$$F_y(y^+) \leq F_y(x),\; \forall x\in  \R^n.$$
By definition of subdifferential at a point, we have that the previous inequality implies
\begin{equation}\label{equ:zero:in:partial:F}
    0 \in \partial F_y(y^+).
\end{equation}
We have that $\Psi_{\cX}$ is a proper closed function and $\ri{\dom{\Psi_{\cX}}}\neq \emptyset$ (see the first two claims of Property \ref{prop:characterization:Psi:cX}). The domain of the quadratic function $h_y:\R^n\to \R$ is $\R^n$. Since $h_y$ is a continuous real value function in $\R^n$, it is also closed (see Proposition 1.1.3 in \cite{Bertsekas:09}). We have that
\begin{align*}
    \ri{\dom{\Psi_{\cX}}}\bigcap \ri{\dom{h_y}} &= \ri{\dom{\Psi_{\cX}}} \bigcap \R^n \\
    &= \ri{\dom{\Psi_{\cX}}}\neq \emptyset.
\end{align*}
Since $F_y=\Psi_{\cX}+h_y$ is equal to the sum of two closed convex functions and $$\ri{\dom{\Psi_{\cX}}}\bigcap \ri{\dom{h_y}}\neq \emptyset,$$ we have
$\partial F_y(y^+) = \partial \Psi_{\cX}(y^+) + \partial h_y(y^+)$ (see Proposition 5.4.6 in \cite{Bertsekas:09}).
The subdifferential of the differentiable function $h_y$ at $y^+$ is $\dhh_y(y^+) = \dhh(y)+R(y^+-y)$. Thus, we obtain from (\ref{equ:zero:in:partial:F})
\begin{align*}
	 0 \in \partial F_y(y^+) &= \partial \Psi_{\cX}(y^+) + \partial h_y(y^+) \\
	 &= \partial \Psi_{\cX}(y^+) + \dhh(y)+R(y^+-y).
\end{align*}
Since $g(y)$ is defined as $R(y-y^+)$ we obtain
$$ g(y)-\dhh(y) \in \partial \Psi_{\cX}(y^+).$$
By definition of $\partial \Psi_{\cX}(\cdot)$ we have
$$ \Psi_{\cX}(x) \geq \Psi_{\cX}(y^+) + \sp{g(y)-\dhh(y)}{x-y^+}, \;\; \forall x\in \R^n.$$
Obviously, since $\cX\subseteq \R^n$, this implies
$$ \Psi_{\cX}(x) \geq \Psi_{\cX}(y^+) + \sp{g(y)-\dhh(y)}{x-y^+}, \;\; \forall x\in \cX.$$
Since $y^+\in \cX$ and $\Psi_{\cX}=\Psi$ for every $x\in \cX$, we obtain
\begin{equation}\label{ineq:sub:gradient}
    \Psi(x) \geq \Psi(y^+) + \sp{g(y)-\dhh(y)}{x-y^+}, \;\forall x\in \cX.
\end{equation}
The convexity of $h(\cdot)$ implies
$$ h(x) \geq h(y) + \sp{\dhh(y)}{x-y}, \;\; \forall x\in \cX.$$
Adding this inequality to (\ref{ineq:sub:gradient}) yields
\begin{eqnarray}
    f(x) &\mkern-14mu = &\mkern-14mu \Psi(x) + h(x) \nonumber\\
         &\mkern-14mu \geq &\mkern-14mu \Psi(y^+) + \sp{g(y)-\dhh(y)}{x-y^+} \nonumber \\
         &&\mkern-14mu +  h(y) + \sp{\dhh(y)}{x-y} \nonumber\\
         &\mkern-14mu = &\mkern-14mu \Psi(y^+) +\sp{g(y)}{x-y^+} \nonumber \\
         &&\mkern-14mu +  h(y) + \sp{\dhh(y)}{y^+-y}, \; \forall x\in \cX. \label{ineq:f:one}
\end{eqnarray}
From Assumption \ref{assum:conv:smooth} we have
\begin{align*}
	h(y) &\geq h(y^+) -\sp{\dhh(y)}{y^+-y} -\frac{1}{2} \|y^+-y\|^2_R \\
	& = h(y^+) -\sp{\dhh(y)}{y^+-y} -\frac{1}{2} \|R^{-1}g(y)\|^2_R \\
	& = h(y^+) -\sp{\dhh(y)}{y^+-y} -\frac{1}{2} \|g(y)\|_*^2 .
\end{align*}
Adding this inequality to  (\ref{ineq:f:one}) yields
\begin{align*}
	f(x) & \geq \Psi(y^+) + h(y^+) +\sp{g(y)}{x-y^+} - \frac{1}{2} \|g(y)\|_*^2\\
	& = f(y^+) +\sp{g(y)}{x-y^+} - \frac{1}{2} \|g(y)\|_*^2, \; \forall x\in \cX.
\end{align*}
From this inequality we have
$$f(y^+) -f(x) \leq \sp{g(y)}{y^+-x} + \frac{1}{2}\|g(y)\|_*^2, \;\forall x\in \cX.$$
This proves (\ref{ineq:cg:a}). We now prove (\ref{ineq:cg:b}) and (\ref{ineq:cg:c}) by means of simple algebraic manipulations.
\small
\begin{eqnarray}
    f(y^+) -f(x) \mkern-34mu &&\leq \sp{g(y)}{y^+-x} + \frac{1}{2}\|g(y)\|_*^2 \nonumber \\
                 &&= \sp{g(y)}{y-x+y^+-y} + \frac{1}{2}\|g(y)\|_*^2 \nonumber \\
                 &&= \sp{g(y)}{y-x} + \sp{g(y)}{y^+-y} + \frac{1}{2}\|g(y)\|_*^2 \nonumber \\
                 &&= \sp{g(y)}{y-x} + \sp{g(y)}{-R^{-1}g(y)} + \frac{1}{2}\|g(y)\|_*^2 \nonumber \\
                 &&= \sp{g(y)}{y-x}- \|g(y)\|_*^2 + \frac{1}{2}\|g(y)\|_*^2 \nonumber \\
                 &&= \sp{g(y)}{y-x}- \frac{1}{2}\|g(y)\|_*^2, \; \forall x\in \cX. \label{ineq:proof:ff}
\end{eqnarray} \normalsize
This proves (\ref{ineq:cg:b}). From this inequality, and the definition of $g(y)$, we obtain
\begin{align*}
	f(y^+) -f(x) & \leq \sp{R(y-y^+)}{y-x} - \frac{1}{2} \| R(y-y^+)\|_*^2 \\
	& = -\sp{R(y-y^+)}{x-y} - \frac{1}{2} \| y-y^+\|_R^2 \\
	& = - \frac{1}{2}\| y-y^+ + x-y\|^2_R + \frac{1}{2}\|x-y\|_R^2 \\
	& = - \frac{1}{2}\| y^+ -x \|^2_R + \frac{1}{2}\|y-x\|_R^2, \; \forall x\in \cX.
\end{align*}
This proves (\ref{ineq:cg:c}). Suppose now that $y\in \cX$. Particularizing inequality (\ref{ineq:proof:ff}) to $x=y$ yields
$$ \frac{1}{2} \| g(y)\|_*^2 \leq f(y)-f(y^+), \;\; \forall y\in \cX.$$
The inequality $f(y)-f(y^+) \leq f(y) -f^*$ trivially follows from $f^*\leq f(y^+)$.
\QED

\subsection{Characterization of optimality}\label{appen:optimaliticharacterization}

The following property serves to characterize the optimality of a given point $y\in \R^n$.

\begin{property}\label{prop:charac:optimality}
	Suppose that Assumption \ref{assum:conv:smooth} holds. Then $y\in \R^n$ belongs to the optimal set 
	$$\Omega=\set{x}{x\in \cX, f(x)=f^*}$$ if and only if $g(y)=0$.
\end{property}  

\proof
We first show that $g(y)=0$ implies $y\in \Omega$. Since $R\succ 0$, we infer from equality $g(y)=R(y-y^+)$ that ${g(y)=0}$ is equivalent to $y=y^+$. Suppose that ${x^*\in \Omega \subseteq \cX}$. Then, we obtain from $g(y)=0$, $y=y^+\in \cX$, and the first claim of Property \ref{prop:proj:gradient}, the following inequality
\begin{align*}
	f(x^*) & \geq f(y^+) -\sp{g(y)}{y^+-x^*}-\frac{1}{2}\|g(y)\|_*^2 \\
	& = f(y^+) = f(y). 
\end{align*}
That is, $f^*=f(x^*)\geq f(y)$. Since $y = y^+ \in \cX$, this is possible only if $y$ is also optimal ($f(y)=f^*$).
This proves that ${g(y)=0}$ implies $y\in \Omega$. We now prove that $y\in \Omega$ implies $g(y)=0$. Suppose that $y\in \Omega$. Then, $f(y)=f^*$ and we obtain from the second claim of  Property \ref{prop:proj:gradient} 
$$ \frac{1}{2}\|g(y)\|_*^2 \leq f(y)-f^*=0.$$
This implies $g(y)=0$.  
\QED

\subsection{Convergence of non restart FISTA \label{app:conver:Fista}}

\begin{property}\label{prop:convergence}
Suppose that Assumption \ref{assum:conv:smooth} holds. Then, the sequences $\{x_k\}$ and $\{y_k\}$ generated by Algorithm \ref{alg:FISTA} (FISTA) satisfy
\blista
\item $  f(x_{k})-f^* \leq \fracg{2\|x_0-\bx_0\|_R^2}{(k+1)^2}$, for all $k\geq 1$,
\item $ \|g(y_k)\|_* \leq \fracg{4\| x_0-\bx_0\|_R}{k+2} $, for all $k\geq 0$,
\elista
where $\bar{x}_0$ represents the point in the optimal set $\Omega$ closest to the initial condition $x_0$ of the algorithm.
\end{property}

\proof 

{\bf First claim:}

We denote $g_k\doteq g(y_k)$, $\forall k\geq 0$. Additionally, we recall that ${\|\cdot\|_* \doteq \|\cdot\|_{R^{-1}}}$.

From step 4 of FISTA algorithm we have
\begin{equation}\label{equ:recu:x:comp}
    x_k =y_{k-1}^+, \; \forall k\geq 1.
\end{equation}
This implies that
$$ g_k = R(y_k-y_k^+) = R(y_k-x_{k+1}),\; \forall k\geq 0.$$
Particularizing inequality (\ref{ineq:cg:c}) of the first claim of Property \ref{prop:proj:in:appendix} to $y=y_0\in\R^n$, and $x=\bx_0\in \Omega \subseteq \cX$, we obtain
$$ f(y_0^+)-f(\bx_0) \leq -\frac{1}{2} \| y_0^+ -\bx_0\|^2_R + \frac{1}{2} \|y_0-\bx_0\|^2_R.$$
By construction we have that $x_0=y_0$ and $x_1=y_0^+$. Furthermore, by definition of $\bx_0$, we have $f(\bx_0)=f^*$. Therefore we can rewrite previous inequality as
\begin{align} 
    f(x_1)-f^* & \leq -\frac{1}{2} \|x_1-\bx_0\|^2_R + \frac{1}{2} \|x_0-\bx_0\|^2_R  \label{ineq:xone}\\
    & \leq \frac{1}{2} \|x_0-\bx_0\|^2_R. \nonumber 
\end{align}
This proves the  claim of the property for $k=1$. We now proceed to prove the claim for $k\geq 2$. From equality (\ref{equ:recu:x:comp}) we have
$$ x_{k+1} = y_k^+, \; \forall k\geq 1.$$
Therefore, from inequality (\ref{ineq:cg:b}) of Property \ref{prop:proj:in:appendix} we obtain that for every $x\in \cX$ and every $k\geq 1$
\begin{eqnarray*}
	f(x) & \geq & f(x_{k+1})+\frac{1}{2}\|g_k\|^2_*- \sp{g_k}{y_k-x}.
\end{eqnarray*}
We notice that, by construction, $x_k\in \cX$, $k\geq 1$. Particularizing at $x_k$ and $\bx_0$, we obtain from last inequality
\begin{subequations}
\begin{flalign}
    f(x_k) & \geq  f(x_{k+1})+\frac{1}{2}\|g_k\|^2_* - \sp{g_k}{y_k-x_k}, \; \forall k\geq 1, \label{ineq:no:increm:a}
    \\ f(\bx_0) & \geq  f(x_{k+1})+\frac{1}{2}\|g_k\|^2_* - \sp{g_k}{y_k-\bx_0}, \; \forall k\geq 1. \label{ineq:no:increm:b}
\end{flalign}
\end{subequations}
In order to write down the proof in a compact way, we introduce the following incremental notation, valid for all $k\geq 0$,
\begin{align*}
    \delta f_k & \doteq f(x_k) - f^*, \\
	\delta x_k & \doteq x_k -\bx_0, \\
	\delta y_k & \doteq y_k -\bx_0, .
\end{align*}
Inequalities (\ref{ineq:no:increm:a}) and (\ref{ineq:no:increm:b}) in an incremental notation, are
\begin{subequations}
\begin{align}
    \delta f_k -\delta f_{k+1} & \geq  \frac{1}{2}\|g_k\|^2_* - \sp{g_k}{\delta y_k-\delta x_k},\; \forall k\geq 1, \label{ineq:f:x:plus} \\
    - \delta f_{k+1} & \geq  \frac{1}{2}\|g_k\|^2_* - \sp{g_k}{\delta y_k}, \; \forall k\geq 1. \label{ineq:f:x:os}
\end{align}
\end{subequations}
We introduce now the auxiliary variable $\Gamma_k$, defined as
$$\Gamma_k\doteq  t_{k-1}^2\delta f_k-t_k^2 \delta f_{k+1}, \; \forall k\geq 1.$$
From Property \ref{prop:tk} in appendix \ref{appen:t} we have
$$t_{k-1}^2=t_k^2-t_k, \; \forall k\geq 1.$$
We now use this identity to obtain
\begin{align}
    \Gamma_k & = (t_k^2-t_k)\delta f_k-t_k^2 \delta f_{k+1} \nonumber\\
    & = (t_{k}^2-t_{k})(\delta f_k-\delta f_{k+1}) - t_k \delta f_{k+1}, \; \forall k\geq 1. \label{equ:aux:Gam}
\end{align}
In view of Property \ref{prop:tk}, $t_k\geq 1$, $\forall k\geq 0$. This implies that we can replace, in inequality (\ref{equ:aux:Gam}), $\delta f_k-\delta f_{k+1}$ and $-\delta f_{k+1}$ by the lower bounds given by  inequalities (\ref{ineq:f:x:plus}) and (\ref{ineq:f:x:os}). In this way we obtain
\begin{eqnarray}
    {}\mkern-52mu \Gamma_k &\mkern-14mu \geq &\mkern-14mu (t_{k}^2-t_{k}) \left(\frac{1}{2}\|g_k\|^2_* - \sp{g_k}{\delta y_k-\delta x_k} \right) \nonumber \\
    &&\mkern-14mu + t_k \left(\frac{1}{2}\|g_k\|^2_* - \sp{g_k}{\delta y_k} \right) \nonumber\\
    &\mkern-14mu = &\mkern-14mu \frac{t_k^2 }{2}\|g_k\|_*^2 - \sp{g_k}{t_k^2(\delta y_k -\delta x_k)+t_k\delta x_k}, \forall k \geq 1. \label{ineq:Gamma:k}
\end{eqnarray}
From step 6 of the algorithm we have for all $k\geq 1$ that ${y_k = x_k + \fracg{t_{k-1} - 1}{t_k} (x_k - x_{k-1})}$.
This can be rewritten in incremental notation as
\begin{equation}\label{yk:in:incremental}
    \delta y_k -\delta x_k = \fracg{t_{k-1} - 1}{t_k} (\delta x_k - \delta x_{k-1}), \; \forall k\geq 1.
\end{equation}
We now define, for every $k\geq 1$
\begin{equation}
    s_k \doteq \delta x_{k-1} + t_{k-1}(\delta x_k -\delta x_{k-1}). \label{equ:sk:one}
\end{equation}
From the definition of $s_k$ and (\ref{yk:in:incremental}) we obtain
\begin{align}
    s_k -\delta x_k & = \delta x_{k-1} +t_{k-1}(\delta x_k-\delta x_{k-1}) -\delta x_{k}  \nonumber \\
    & = (t_{k-1} -1) (\delta x_k-\delta x_{k-1}) \nonumber \\
    & = t_k( \delta y_k - \delta x_k), \; \forall k\geq 1. \label{equ:sk:two}
\end{align}
From (\ref{ineq:Gamma:k}) and (\ref{equ:sk:two})  we obtain
\begin{align}
    \Gamma_k & \geq \frac{1}{2}\|t_k g_k\|_*^2 - \sp{g_k}{t_k(s_k-\delta x_k)  + t_k\delta x_k} \nonumber\\
    & = \frac{1}{2}\|t_k g_k\|_*^2 -\sp{ t_k g_k}{s_k}, \;\forall k\geq 1.  \label{equ:Gammakupps}
\end{align}
Using (\ref{equ:sk:one}) and (\ref{equ:sk:two}) we now show that $g_k$ can be written in terms of $s_k$ and $s_{k+1}$.
\begin{align}
t_k g_k  & = t_k R(y_k-x_{k+1}) = t_k R(\delta y_k -\delta x_{k+1})  \nonumber\\
         & = t_k R(\delta y_k - \delta x_{k} + \delta x_{k} - \delta x_{k+1}) \nonumber \\
         & = R(s_k-\delta x_k + t_k(\delta x_k -\delta x_{k+1} )) \nonumber \\
         & = R(s_k-s_{k+1}), \; \forall k\geq 1. \label{equ:gk:tk:equ}
\end{align}
With this expression for $t_k g_k$ we obtain from (\ref{equ:Gammakupps})
\begin{align*}
	\Gamma_k & \geq \frac{1}{2}\|R (s_k-s_{k+1})\|_*^2 -\sp{R(s_k-s_{k+1})}{s_k} \\
	& = \frac{1}{2}\|s_{k+1}-s_k\|^2_R +\sp{R(s_{k+1}-s_k)}{s_k} \\
	& = \frac{1}{2}\|(s_{k+1}-s_k)+s_k\|^2_R - \frac{1}{2}\|s_k\|^2_R \\
	& = \frac{1}{2} \|s_{k+1}\|^2_R - \frac{1}{2}\|s_{k}\|^2_R, \; \forall k\geq 1.
\end{align*}
Thus, for every $k\geq 1$,
$$ \Gamma_k = t_{k-1}^2 \delta f_k - t_k^2 \delta f_{k+1} \geq \frac{1}{2} \|s_{k+1}\|^2_R - \frac{1}{2}\|s_{k}\|^2_R. $$
Equivalently
$$ t_k^2 \delta f_{k+1}+\fracg{1}{2}\|s_{k+1}\|^2_R \leq t_{k-1}^2\delta f_k +\fracg{1}{2}\|s_k\|^2_R, \; \forall k\geq 1.$$
Since this inequality holds for every $k\geq 1$ we can apply it in a recursive way to obtain
\begin{equation*}
    \begin{aligned}
	t_k^2 \delta f_{k+1} + \frac{1}{2}\|s_{k+1}\|^2_R & \leq t_0^2 \delta f_1 + \frac{1}{2} \| s_1\|^2_R \\
	& =  \delta f_1+  \frac{1}{2} \| \delta x_0 + t_0(\delta x_1-\delta x_0)\|^2_R \\
	& =  \delta f_1 + \frac{1}{2} \| x_1 - \bx_0\|^2_R,\;  \forall k\geq 1.
\end{aligned}
\end{equation*}
\vspace{-0.5em}
From (\ref{ineq:xone}) we have $${f(x_1)-f^* + \frac{1}{2} \| x_1 - \bx_0\|^2_R \leq \frac{1}{2}\|x_0-\bx_0\|^2_R}.$$
\vspace{-0.5em}
Thus,
\begin{equation}\label{ineq:convergence}
t_k^2 \delta f_{k+1} + \frac{1}{2}\|s_{k+1}\|^2_R  \leq \frac{1}{2}\| x_0-\bx_0\|^2_R, \; \forall k\geq 1.
\end{equation}
Therefore,
\vspace{-0.5em}
$$t_k^2 (f(x_{k+1})-f^*) + \frac{1}{2} \|s_{k+1}\|^2_R \leq \frac{1}{2} \|x_0-\bx_0\|^2_R, \; \forall k\geq 1.$$
From this inequality, and taking now into account that ${t_k \geq \fracg{k+2}{2}}$ for all $k\geq 0$ (second claim of Property \ref{prop:tk}), we conclude
$$ f(x_{k+1})-f^* \leq \frac{\|x_0-\bx_0\|^2_R}{2 t_k^2} \leq \frac{2 \|x_0-\bx_0\|^2_R}{(k+2)^2}, \; \forall k\geq 1.$$
That is,
\vspace{-0.5em}
$$ f(x_{k})-f^* \leq \frac{2 \|x_0-\bx_0\|^2_R}{(k+1)^2}, \;\; \forall k\geq 2.$$

\QED

{\bf Second claim:}

We first prove the claim for $k=0$.
\begin{align*}
	\|g(y_0)\|_* & = \| R(y_0-y_0^+)\|_* = \| y_0-y_0^+\|_R \\
                 & = \| x_0 -x_1\|_R = \| x_0-\bx_0 +\bx_0 -x_1 \|_R\\
                 & \leq \| x_0 -\bx_0\|_R + \| x_1-\bx_0\|_R.
\end{align*}
From (\ref{ineq:xone}) we derive
\begin{equation} \label{ineq:xone:closer:xzero}
    \|x_1-\bx_0\|_R \leq \|x_0- \bx_0\|_R.
\end{equation}
Thus,
$$ \|g(y_0)\|_* \leq \| x_0 -\bx_0\|_R + \| x_1-\bx_0\|_R \leq 2\|x_0-\bx_0\|_R.$$
We now prove the claim for $k>0$. From (\ref{ineq:convergence}) we also have
\begin{equation}\label{ineq:norm:s:plus}
    \|s_{k+1}\|_R \leq  \|x_0-\bx_0\|_R, \; \forall k\geq 1.
\end{equation}
We also have that
\begin{equation} \label{equ:sone:dxzero}
    s_1=\delta x_0 + t_0(\delta x_1-\delta x_0) = x_1-\bx_0.\end{equation}
From (\ref{ineq:xone:closer:xzero}) we derive
${\|s_1\|_R=\| x_1-\bx_0\|_R \leq \|x_0-\bx_0\|_R}$. From this and (\ref{ineq:norm:s:plus}) we obtain
\begin{equation}\label{ineq:norm:sk}
    \|s_{k}\|_R \leq  \|x_0-\bx_0\|_R, \; \forall k\geq 1.
\end{equation}
From here we derive, for every $k\geq 1$,
\begin{flalign*}
	\| s_{k+1}-s_k \|_R & \leq \| s_{k+1} \|_R + \| s_k \|_R \\
	& \leq \| x_0 -\bx_0 \|_R + \| x_0-\bx_0\|_R = 2 \| x_0-\bx_0\|_R.
\end{flalign*}
From (\ref{equ:gk:tk:equ}) we have
$$ g_k = \frac{1}{t_k} R(s_k-s_{k+1}), \forall k\geq 1.$$
Therefore, for every $k\geq 1$
\begin{align*}
	\| g_k\|_* & = \frac{1}{t_k} \| s_k-s_{k+1}\|_R \\
	& \leq \frac{2}{t_k}\|x_0-\bx_0\|_R\\
	& \leq \frac{4}{k+2}\|x_0-\bx_0\|_R.
\end{align*}
We notice that the last inequality is due to the second claim of Property \ref{prop:tk}. 
This proves the second claim of the property.
\QED

\subsection{Properties of the sequence $\{t_k\}$} \label{appen:t}
\begin{property}\label{prop:tk}
	Let us suppose that $t_0=1$ and that
	$$t_k \doteq \frac{1}{2} \left( 1+\sqrt{1+4t_{k-1}^2}\,\,\right), \;\; \forall k\geq 1.$$
	Then
	\blista
	\item $t_{k-1}^2=t_k^2-t_k$, for all $k\geq 1$.
	\item $ t_k \geq \fracg{k+2}{2}\geq 1$, for all $k\geq 0$.
	\elista
\end{property}

\proof \mbox{}

\blista
\item For every $k\geq 1$, $t_k$ is defined as one of the roots of
$$t_k^2-t_k-t_{k-1}^2=0.$$
Therefore we obtain $t_{k-1}^2=t_k^2-t_k$.
\item The claim is trivially satisfied for $k$ equal to 0. We now show that if the claim is satisfied for $k-1$ then it is also satisfied for $k$.
\begin{align*}
	t_k &= \frac{1}{2} \left( 1+\sqrt{1+4t_{k-1}^2}\,\,\right) \\
        & \geq \frac{1}{2} \left( 1+ \sqrt{4 t_{k-1}^2} \right) = \frac{1}{2} + t_{k-1}.
\end{align*}
Since the claim is assumed to be satisfied for $k-1$ we have $t_{k-1}\geq \frac{k+1}{2}$ and consequently
$$ t_k \geq \frac{1}{2}+\frac{k+1}{2} = \frac{k+2}{2}.$$
\QED
\elista

\subsection{Proof of Property \ref{prop:FISTA:EC}}\label{appen:FISTA:EC}

From equation (\ref{equ:conv:FISTA}) we have
\begin{equation*} f(x_k)-f^* \leq \frac{2}{(k+1)^2}\|x_0-\bx_0\|_{R}^2, \; \forall k\geq 1.
\end{equation*}
Due to Assumption \ref{assump:quadratic:growth} we also have
$$ \frac{\mu}{2} \|x_0-\bx_0\|_R^2 \leq f(x_0)-f^*.$$ Therefore,
\begin{equation}\label{equ:pre:alpha}
 f(x_k)-f^* \leq \frac{4}{\mu(k+1)^2}(f(x_0)-f^*),\; \forall k\geq 1.
 \end{equation}
This proves the first claim. Denote 
$$\alpha_k\doteq\fracg{4}{\mu(k+1)^2}, \;\; \forall k\geq 1.$$
With this notation we rewrite (\ref{equ:pre:alpha}) as 
\begin{equation}\label{ineq:with:alpha:k}
    f(x_k)-f^* \leq \alpha_k (f(x_0)-f^*), \; \forall k\geq 1.
\end{equation}
Suppose now that  $k\geq \left\lfloor \fracg{2}{\sqrt{\mu}}\right\rfloor$. Then, 
\small
\begin{equation*}
    \alpha_k =\frac{4}{\mu (k+1)^2} \leq  \fracg{4}{\mu \left( \left\lfloor \fracg{2}{\sqrt{\mu}}\right\rfloor +1\right)^2} < \frac{4}{\mu\left(\fracg{2}{\sqrt{\mu}}\right)^2} = 1.
\end{equation*}
\normalsize
Therefore, 
\begin{equation}\label{equ:alpha:in:zeroone}
    \alpha_k\in(0,1),\;  \forall k\geq \left\lfloor \fracg{2}{\sqrt{\mu}}\right\rfloor.  
\end{equation}
This, along with inequality (\ref{ineq:with:alpha:k}), yields
$$ f(x_k)-f^* \leq f(x_0)-f^*, \;  \forall k\geq \left\lfloor \fracg{2}{\sqrt{\mu}}\right\rfloor.$$
Equivalently,
$$ f(x_k) \leq f(x_0), \;  \forall k\geq \left\lfloor \fracg{2}{\sqrt{\mu}}\right\rfloor. $$
This proves the second claim of the property. 
In view of inequality (\ref{ineq:with:alpha:k}) we have 
\begin{align*}
f(x_k)-f^* & \leq \alpha_k (f(x_0)-f^*) \\
           & = \alpha_k (f(x_0)-f(x_k)+f(x_k)-f^*) \\
           & = \alpha_k (f(x_0)-f(x_k)) + \alpha_k(f(x_k)-f^*).
\end{align*}
Therefore,
\begin{equation} \label{ineq:alpha:k:f} 
    (1-\alpha_k) (f(x_k)-f^*) \leq \alpha_k(f(x_0)-f(x_k)).
\end{equation}
Suppose now that $k\geq \left\lfloor \frac{2\sqrt{e+1}}{\sqrt{\mu}}\right\rfloor$. This implies $k\geq \left\lfloor \frac{2}{\sqrt{\mu}}\right\rfloor$ and consequently $1-\alpha_k>0$ (see (\ref{equ:alpha:in:zeroone})).
Dividing both terms of inequality (\ref{ineq:alpha:k:f}) by $1-\alpha_k$, we get
\begin{align*}
    f(x_k)-f^* & \leq \frac{\alpha_k}{1-\alpha_k} (f(x_0)-f(x_k))\\
    & = \frac{\frac{4}{\mu(k+1)^2}}{1-\frac{4}{\mu(k+1)^2}} (f(x_0)-f(x_k))\\
    & = \frac{4(f(x_0)-f(x_k))}{\mu(k+1)^2 -4} \\ 
    & \leq \fracg{4(f(x_0)-f(x_k))}{\mu(\left\lfloor \frac{2\sqrt{e+1}}{\sqrt{\mu}} \right\rfloor+1)^2-4} \\
    & \leq \fracg{4(f(x_0)-f(x_k))}{\mu(\frac{2\sqrt{e+1}}{\sqrt{\mu}})^2-4} \\
    & = \fracg{4(f(x_0)-f(x_k))}{4(e+1)-4} = \frac{f(x_0)-f(x_k)}{e}.
\end{align*}
\QED
 
 \subsection{Proof of Property \ref{prop:L:Convergence:FISTA}}\label{appen:L:Convergence:FISTA}

By construction, $r_{j-1}\in \cX$, for all $j \geq 1$. Therefore,  
we have from the second claim of Property \ref{prop:proj:gradient}, that 
\begin{equation}\label{ineq:grad:r:r:minus}
    \frac{1}{2} \|g(r_{j-1})\|_*^2 \leq f(r_{j-1})-f(r_{j-1}^+), \; \forall j\geq 1.
\end{equation}  
We also notice that  $r_j$ is computed invoking FISTA algorithm using $r_{j-1}$ as initial condition ($z=r_{j-1}$). That is, 
$$[r_j,n_j]=FISTA(r_{j-1},n_{j-1},E_c^l).$$ 
Since the output value $f(r_j)$ is forced to be no larger than the one corresponding to ${x_0=z^+=r_{j-1}^+}$, we have ${f(r_{j})\leq f(r_{j-1}^+)}$. Therefore, we obtain from inequality (\ref{ineq:grad:r:r:minus}) that 
\begin{align*}
    \frac{1}{2} \|g(r_{j-1})\|_*^2 & \leq f(r_{j-1})-f(r_{j-1}^+)\\
                                   & \leq f(r_{j-1})-f(r_{j}).
\end{align*}
This proves the first claim of the property.  
We now show that if $n_{j-1} \leq \frac{4\sqrt{e+1}}{\sqrt{\mu}}$, then the value $n_j$ obtained from 
$$ [r_j,n_j]=FISTA(r_{j-1},n_{j-1},E_c^l),$$ 
also satisfies 
\begin{equation} \label{eq:n_geq_4e}
    n_{j} \leq \frac{4\sqrt{e+1}}{\sqrt{\mu}}.
\end{equation}
Denote 
$$\bm=\left\lfloor \frac{2\sqrt{e+1}}{\sqrt{\mu}} \right\rfloor.$$
Since  $\bm \geq \left\lfloor \frac{2 \sqrt{e+1}}{\sqrt{\mu}} \right\rfloor$,  we infer, from the third claim of Property \ref{prop:FISTA:EC}, that
$$ f(x_{\bm})-f^* \leq \frac{f(x_0)-f(x_{\bm})}{e}. $$ 
From this inequality, we obtain 
$$ f(x_{\bm})-f(x_k) \leq f(x_{\bm})-f^* \leq \frac{f(x_0)-f(x_{\bm})}{e}.$$
Therefore, the first exit condition is satisfied for $m=\bm$.
Since $m=\lfloor \frac{k}{2} \rfloor+1$ we have $m \geq \frac{k}{2}$. 
This means that for $m=\bm$, the corresponding value for $k$ is no larger than 
$$2\bm = 2 \left\lfloor \frac{2\sqrt{e+1}}{\sqrt{\mu}} \right\rfloor \leq \frac{4(\sqrt{e+1})}{\sqrt{\mu}}.$$
We also notice that, in view of the second claim of Property \ref{prop:FISTA:EC}, the additional exit condition $f(x_k)\leq f(x_0)$ is satisfied for every 
$$k\geq \left \lfloor \frac{2}{\sqrt{\mu}} \right\rfloor.$$ 
Therefore, $n_{j-1}\leq \frac{4\sqrt{e+1}}{\sqrt{\mu}} $ implies that 
$n_j$, obtained from $[r_j,n_j]=FISTA(r_{j-1},n_{j-1},E_c^l)$, 
also satisfies (\ref{eq:n_geq_4e}).
We now prove, by reduction to the absurd, that $n_j$ cannot be larger than 
$\frac{4\sqrt{e+1}}{\sqrt{\mu}}$. Suppose that 
\begin{equation}\label{ineq:lower:nj}
    n_j > \frac{4\sqrt{e+1}}{\sqrt{\mu}}.
\end{equation}
Because of the previous discussion, the previous inequality could be forced only by the doubling step ${n_j = 2n_{j-1}}$ of the algorithm. 
That is, inequality (\ref{ineq:lower:nj}) is possible only if there is $s$ such that $n_{s-1} >\frac{2\sqrt{e+1}}{\sqrt{\mu}}$ and 
$$ f(r_{s-1})-f(r_s) > \frac{f(r_{s-2})-f(r_{s-1})}{e}.$$
Since 
$$ [r_{s-1},n_{s-1}]=FISTA(r_{s-2},n_{s-2},E_c^l),$$
we have that $r_{s-1}$ is obtained from $r_{s-2}$ applying 
$$n_{s-1} >\frac{2\sqrt{e+1}}{\sqrt{\mu}}$$
iterations of FISTA algorithm. However, we have from the third claim of Property \ref{prop:FISTA:EC} that this number of iterations implies
$$ f(r_{s-1})-f(r_s) \leq f(r_{s-1})-f^* \leq 
\frac{f(r_{s-2}^+)-f(r_{s-1})}{e}.$$
From the second claim of Property \ref{prop:proj:gradient} we also have ${f(r_{s-2}^+)\leq f(r_{s-2})}$. Thus,  
$$ f(r_{s-1})-f(r_s) \leq \frac{f(r_{s-2})-f(r_{s-1})}{e}.$$
That is, there is no doubling step if $n_{s-1}\geq \frac{2\sqrt{e+1}}{\sqrt{\mu}}$. This proves the second claim of the property.

We now show that there is a doubling step at least every 
$$T \doteq \left \lceil\ln\left(1+\frac{2(f(r_0)-f^*)}{\epsilon^2}\right)\right\rceil$$
steps of the algorithm. Suppose that there is no doubling step from iteration $j=s+1$ to $j=s+T$, where $s\geq 1$. That is, 
$$
f(r_{j-1})-f(r_j) \leq \frac{f(r_{j-2}) -f(r_{j-1})}{e}, \; \forall j \in [s+1,s+T].
$$
From this, and the first claim of the property, we obtain the following sequence of inequalities   
\small
\begin{equation*}
    \begin{aligned}
&\frac{1}{2} \|g(r_{s+T-1})\|_*^2 \leq f(r_{s+T-1})-f(r_{s+T})   \\
& \leq \frac{f(r_{s+T-2})-f(r_{s+T-1})}{e} \leq \left( \frac{1}{e} \right)^{T}(f(r_{s-1})-f(r_{s}))\\
& \leq \left( \frac{1}{e} \right)^{T} (f(r_{s-1})-f^*) \leq \left( \frac{1}{e} \right)^{T} (f(r_0)-f^*) \\
& = \left( \frac{1}{e} \right)^{\left\lceil \ln \left(1+\frac{2(f(r_0)-f^*)}{\epsilon^2}\right)\right\rceil} (f(r_0)-f^*) \\
& \leq \left( \frac{1}{e} \right)^{\ln \left(1+\frac{2(f(r_0)-f^*)}{\epsilon^2}\right)} (f(r_0)-f^*) \\
& = \left( \fracg{1}{1+\frac{2(f(r_0)-f^*)}{\epsilon^2}}\right) (f(r_0)-f^*) \leq \frac{\epsilon^2}{2}.
\end{aligned}
\end{equation*}
\normalsize
We conclude that $T$ consecutive iterations without doubling step implies that the exit condition is satisfied ($\|g(r_{s+T-1})\|_*\leq \epsilon$). We conclude that there must be at least one doubling step every $T$ iterations. This implies that there exist $j\in[s+1,s+T]$ such that 
$$ f(r_{j-1})-f(r_j) > \frac{f(r_{j-2}) -f(r_{j-1})}{e}.$$
Therefore, $n_{j} = 2n_{j-1}$. Moreover, since $\{n_j\}$ is a non decreasing sequence, we get ${n_{s+T} \geq n_j = 2n_{j-1} \geq 2n_{s}}$, ${\forall s\geq 1}$.
That is,
\begin{equation}\label{ineq:half:n:T}
n_s\leq \frac{n_{s+T}}{2}, \; \forall s\geq 1.
\end{equation}
Suppose that $j$ is rewritten as $j=m+nT$, where $0\leq m<T$ and $n\geq 0$. From the non decreasing nature of $\{n_j\}$,
\small
\begin{equation}
    \begin{aligned}
     \Sum{i=0}{j}& n_i = \Sum{i=0}{m+nT} n_i = \Sum{i=0}{m} n_i + \Sum{\ell=0}{n-1} \Sum{i=1}{T} n_{m+i+\ell T} \label{ineq:sum:nj}
\\
 \leq& T n_m +  T \Sum{\ell=1}{n} n_{m+\ell T} = T\Sum{\ell=0}{n} n_{m+\ell T} = T\Sum{\ell=0}{n} n_{j-\ell T}.
\end{aligned} 
\end{equation}
\normalsize%
Also, from inequality (\ref{ineq:half:n:T}), we have $n_{j-T} \leq \frac{n_j}{2}$. Using this inequality in a recursive manner we obtain
$$ n_{j-\ell T} \leq \left(\frac{1}{2}\right)^{\ell} n_j, \; \ell =0,\ldots,n.$$
This, allows us to infer from (\ref{ineq:sum:nj}) that
\begin{equation*}
    \Sum{i=0}{j} n_i \leq T \Sum{\ell=0}{n} \left(\frac{1}{2}\right)^\ell n_{j} \leq T \Sum{\ell=0}{\infty} \left(\frac{1}{2}\right)^\ell n_{j} = 2Tn_j.
\end{equation*}
The last claim of the property follows directly from this one and the bound $n_j\leq \frac{4\sqrt{e+1}}{\sqrt{\mu}}$ of the second claim. 
That is, if $j$ denotes the first index for which $\|g(r_j)\|_*\leq \epsilon$, we get that the number of total iterations is bounded by 
\begin{align*}
    \Sum{i=0}{j} n_i & \leq 2Tn_j  \leq  \frac{8T\sqrt{e+1}}{\sqrt{\mu}} \leq  \frac{16T}{\sqrt{\mu}} \\
    & = \frac{16}{\sqrt{\mu}} \left\lceil \ln \left(1+\frac{2(f(r_0)-f^*)}{\epsilon^2} \right)\right\rceil.
\end{align*}
\QED

\end{appendix}

\end{document}